\documentclass[reqno,A4paper]{amsart}
\usepackage{amsmath,amssymb,amscd,amsthm,verbatim,alltt,amsfonts,array}
\usepackage{mathrsfs}
\usepackage[english]{babel}
\usepackage{latexsym}
\usepackage{amssymb}
\usepackage{euscript}
\usepackage{graphicx}
\usepackage{hyperref}
\usepackage{cite}
\usepackage{rotating}
\usepackage{xcolor}
\usepackage{multirow}
\usepackage{longtable}
\usepackage{geometry}
\usepackage{caption}
\usepackage{float}
\newfloat{graphic}{tbp}{lgr}   
\floatname{graphic}{Graphic}
\floatstyle{boxed}
\restylefloat{figure}
\newtheorem{thm}{Theorem}[section]

\newcommand{\be}{\begin{equation}}
\newcommand{\ee}{\end{equation}}
\newcommand{\ben}{\begin{enumerate}}
\newcommand{\een}{\end{enumerate}}
\newcommand{\beq}{\begin{eqnarray}}
\newcommand{\eeq}{\end{eqnarray}}
\newcommand{\beqn}{\begin{eqnarray*}}
\newcommand{\eeqn}{\end{eqnarray*}}
\title[Douglas Curvature in Finsler Geometry]{Douglas Curvature in Finsler Geometry:\\ A Review of the Notion and Its Applications in Geometry and Physics}
\author{Nasrin Sadeghzadeh$^{A}$}
\newcommand{\acr}{\newline\indent}
\address{\llap{a\,}Department of Mathematics,\acr
University of Qom, \acr
Alghadir Bld\acr
Qom\acr
Iran}
\email{nsadeghzadeh@qom.ac.ir}
\author{Meshkat Yavari$^{B}$}
\address{\llap{b\,}Department of Mathematics,\acr
University of Qom, \acr
Alghadir Bld\acr
Qom\acr
Iran}
\email{meshkat.yavari@stu.qom.ac.ir}
\begin{document}
\maketitle
\section*{Abstract}
This review paper provides a comprehensive overview of the current state of research on Douglas curvature in Finsler spaces. It explores the significance, properties, and applications of Douglas curvature, and its role in understanding Finsler geometry. The paper reviews the historical development and significance of Douglas curvature, its characterizations, generalizations, and applications, and presents future research directions.

\textbf{Subjclass} [2010] {53B40; 53C60}

\textbf{Keywords:}{Douglas metrics; Generalized Douglas-Weyl metrics, Douglas metrics of second kind, Conformal transformations of Finsler metrics, Projective transformations.}

\section{Introduction}

Finsler geometry traces its conceptual origins back to Bernhard Riemann's seminal lecture "Ueber die Hypothesen, welche der Geometrie zu Grunde liegen" in 1854 \cite{Rie}. In this work, Riemann laid the foundation for the modern understanding of manifolds, Riemannian structures, and Finsler structures. However, he noted the complexity of the general Finsler case compared to the simpler Riemannian (quadratic) case.

Despite Riemann's foundational work, Finsler geometry remained dormant for decades until it resurfaced in the 1918 thesis of P. Finsler under the guidance of Caratheodory. This resurgence led to the nomenclature of Riemann-Finsler geometry or simply Finsler geometry in honor of its revival by Finsler.

After Finsler's thesis, there was a surge of interest in the field, resulting in the formation of various schools of Finsler geometry and notable contributions from different geometers. This initial phase of development established Finsler geometry as a distinct field within mathematics, often viewed as an extension of Riemannian geometry, \cite{BaoChernShen}, \cite{Miron}, and sparked interest in its utilization in physics as a generalization of pseudo-Riemannian (Lorentzian) geometry \cite{FinsPhysics1}, \cite{FinsPhysics2}.

L. Berwald emerged as a prominent figure during this period, introducing the Berwald connection and a class of spaces bearing his name. These spaces hold significant importance due to their close association with Riemannian spaces, providing a higher level of rigidity that enables the imposition of more stringent conditions. Noteworthy there is a rigidity theorem established by Szab\'{o} which asserts that a Berwald space of dimension 2 is either a Riemannian space or a locally Minkowski space, prompting the need to explore higher dimensions to find examples of Berwald spaces \cite{Szabo}. Moreover, Berwald geometries, being Finsler geometries closely related to (pseudo)-Riemannian geometries, are highly valuable for the application of Finsler geometry in physics \cite{FinsPhysics3}. The proximity of Berwald geometries to general relativity makes them particularly relevant for gravitational theories based on Berwald geometry, which can be viewed as theories closely aligned with general relativity \cite{GR2}. Furthermore, in \cite{GR1}, the exploration of Berwald space-times and their vacuum dynamics is undertaken, with the vacuum dynamics rooted in a Finsler extension of Einstein's equations derived from an action on the unit tangent bundle. The research presented in \cite{FinsPhysics3} focuses on identifying Berwald spaces as a generalization of (pseudo) Riemannian spaces. The paper introduces a simple geometric partial differential equation that determines the conditions under which a given Finsler Lagrangian is of Berwald type. Since Berwald geometries are Finsler geometries that are closely related to (pseudo)-Riemannian geometries, their identification holds great significance for the application of Finsler geometry in physics. In \cite{FinsPhysics3}, the authors illustrate the results with novel examples of $(\alpha, \beta)$-Berwald geometries, which represent Finslerian versions of Kundt (constant scalar invariant) space-times. These findings generalize earlier work by Tavakol and van den Bergh \cite{Tavakol}, as well as the Berwald conditions for Randers and m-Kropina geometries, which are closely related to very special/general relativity geometries.

A Finsler space is said as Berwald type when the Cartan non-linear connection establishes a linear connection on $TM$, or equivalently, an affine connection on the base manifold, where the connection coefficients take the form
\[
N{^i}_j(x,y)=\Gamma{^i}_{jk}(x) y^k,
\]
with $\Gamma: M \rightarrow \mathbb{R}$ denoting a set of smooth functions. The transformation properties of $N{^i}_j$ imply that the functions $\Gamma{^i}_{jk}$ exhibit the appropriate transformation properties to serve as the connection coefficients of a (torsion-free) affine connection on $M$. This affine connection is referred to as the associated affine connection or simply the affine connection within the Berwald space.

Equivalently one can define Berwald geometries by the demand that the geodesic spray $G^i(x,y)=\frac{1}{2}N{^i}_j(x,y) y^j$ be quadratic in $y$. The definition goes back to the original article by Berwald \cite{Ber3}.

Douglas metrics are a generalization of Berwald spaces, a special class of Finsler spaces, introduced by the pioneering work of J. Douglas in the 1920s \cite{Douglas}. It has been studied extensively in the literature, with various aspects of their properties and generalizations being explored.

The notion of Douglas metrics was first introduced by Matsumoto (1997) as a generalization of the concept of Berwald spaces, which are a special class of Finsler spaces \cite{Bacso}.
Douglas metrics form a rich class of Finsler metrics that include locally projectively flat Finsler metrics. They are characterized by the property that the spray coefficients $G^i$ can be expressed in the form
\[
G^i = \frac{1}{2} \Gamma{^i}_{jk}(x) y^j y^k + P(x,y) y^i,
\]
where $\Gamma{^i}_{jk}$ are some local functions and $P$ is a 1-homogeneous function in $y$.

A key result is that a Finsler metric is a Douglas metric if and only if its Douglas tensor, $D_j{^i}_{kl}$ which is defined as follows, vanishes.
\[
D_j{^i}_{kl}=B_j{^i}_{kl}-\frac{2}{n+1}(E_{kl}\delta{^i}_j+E_{jl}\delta{^i}_k+E_{jk}\delta{^i}_l),
\]
where $B_j{^i}_{kl}$ are the components of the Berwald curvature tensor and $E_{jk}=\frac{1}{2}B_j{^m}_{km}$.

Furthermore, as noted in \cite{Sh2}, a spray is projectively affine if $D=0$. Thus, a Finsler metric is called a Douglas metric if its spray is projectively affine. This characterization highlights the essential nature of Douglas metrics within Finsler geometry and connects them to broader geometric properties and implications.

This provides an alternative characterization of Douglas metrics in terms of curvature.
The study of Douglas curvature has been instrumental in unveiling the deeper geometric insights within Finsler spaces, as it provides a lens through which the curvature properties and underlying structures can be analyzed. This review paper aims to offer a comprehensive exploration of the current state of researches on Douglas curvature in Finsler spaces. By synthesizing the findings from seminal works in this field, we seek to elucidate the significance, properties, and applications of Douglas curvature, thereby contributing to the broader understanding of Finsler geometry.
In this review paper, we delve into the intricate world of Douglas curvature in Finsler spaces, aiming to provide a comprehensive overview of the current state of research in this area.
By delving into the existing literature and synthesizing the key findings, this review paper aims to provide a comprehensive understanding of the role of Douglas curvature in Finsler geometry. It will explore the various characterizations, generalizations, and applications of Douglas curvature, ultimately contributing to the broader discourse on the geometric structures and properties of Finsler spaces. This review paper has been structured into several interesting sections.
\ben

\item[]{-} This opening section will provide a comprehensive overview of the historical development and significance of Douglas curvature in the context of Finsler geometry.

\item[]{-} The next section reviews the papers with the subject of Douglas curvature in general form, focusing on the key findings and contributions made in this area.

\item[]{-} The third section  will further explore the specific classes of Finsler metrics with vanishing Douglas curvature, highlighting its role in enhancing our understanding of the curvature properties and underlying structures of Finsler manifolds.

\item[]{-} The fourth section delves into the realm of metric transformations within Douglas spaces. By conducting a comprehensive review of related literature and researches, it analyses the relation between transformations and the fundamental attributes of Douglas metrics. This section is organized into two subsections, "Conformal Transformations and Their Influence on Douglas Spaces" and "Exploring Transformations Beyond Conformal Changes," providing an in-depth analysis of transformations within Douglas spaces.

\item[]{-} Our examination in the fifth part will focus on all papers that elaborate on the concept of Douglas curvature and its generalizations. Throughout this particular portion, we will review the innovative classes of Finsler metrics that have been formulated through these generalizations and extensions of the original Douglas curvature.

\item[]{-} The concluding section of this review paper will present intriguing ideas for future research directions related to Douglas curvature, highlighting potential avenues for further study and application in the field of Finsler geometry.
\een
Each review section, including Sections \ref{Douglas}, \ref{Transformation}, and \ref{Generalization}, is accompanied by tables that summarize the key findings and insights derived from the exploration of Douglas curvature in Finsler metrics and its related generalizations. These tables present the information in a brief and organized manner.

\section{Douglas Curvature in Finsler geometry}\label{Douglas}

This section delves into the concept of Douglas curvature in various papers, starting with the seminal work by B\'{a}cs\'{o} \cite{Bacso}. This paper presents a series of crucial theorems that define Finsler spaces of Douglas type. Specifically, a Finsler space is classified as a Douglas space if the expressions $D^{ij} = G^i y^j - G^j y^i$ are homogeneous polynomials in the coordinate components $y^i$ of degree three. This definition is equivalent to the condition that the Douglas tensor vanishes identically.
The review further explores the relationships between Douglas spaces and other important classes of Finsler metrics. It is shown that if a Finsler space is both a Landsberg space and a Douglas space, then it is necessarily a Berwald space, and conversely, a Berwald space satisfies both the Landsberg and Douglas conditions. Additionally, the paper examines theorems that characterize Douglas Finsler spaces in terms of the projective connection and Wagner connection, revealing that the Douglas curvature is projective invariant. Furthermore, specific conditions are stated for certain classes of Finsler metrics to be of Douglas type. These results provide a comprehensive understanding of the defining properties and characterizations of Finsler metrics that belong to the class of Douglas spaces.

In the study \cite{Bacso2}, a series of theorems shed light on the transformations of Douglas metrics into Berwald metrics. The research reveals that a Landsberg space evolves into a Berwald space when identified as a Douglas space, while conversely, a Douglas space transitions into a Berwald space when displaying Landsberg characteristics. The exploration extends to generalized m-Kropina spaces (distinct from Kropina spaces) and Matsumoto spaces, demonstrating that a Douglas space undergoes a shift to a Berwald space under specific conditions. Furthermore, the analysis delves into the characterization of Finsler spaces concerning metric components, delineating the conditions for a space to be classified as of Douglas type and its subsequent reduction to a Berwald space under particular circumstances. The theorems also investigate the impact of vanishing T-tensors in two-dimensional Douglas spaces, resulting in a conversion to Berwald spaces with a constant main scalar. Lastly, the research examines scenarios where a two-dimensional Douglas space with a cubic metric transforms into either a locally Minkowski space or a Berwald space with distinct properties, offering valuable insights into the dynamics of transformations and classifications in the realm of Finsler geometry.

The papers \cite{DShen} and \cite{DouglasE-cur} provide valuable insights into the relationships between Douglas metrics and other important classes of Finsler metrics, such as Berwald and Landsberg spaces. It highlights the conditions under which Douglas metrics can be transformed or reduced to these related metric types.
\begin{thm}\cite{DShen}

In the case of a non-Riemannian Douglas manifold of dimension $n \geq 3$, the following conditions are equivalent,
\ben
\item [] {(a)} The Finsler metric has isotropic mean Berwald curvature.
\item []{(b)} The Finsler metric has relatively isotropic Landsberg curvature.
\een
\end{thm}

Furthermore, the review covers a theorem from the work by \cite{DouglasE-cur}, which states that for a complete Douglas space with bounded Cartan torsion, if the $\bar{E}$-curvature of the Finsler metric vanishes, then the space reduces to a Berwald metric. In particular, every compact Douglas space with $\bar{E} = 0$ is a Berwald space.
Continuing our exploration of Douglas curvature, we delve into additional research papers focusing on specific classes of Finsler metrics. In the upcoming section, we examine studies that investigate $(\alpha, \beta)$-metrics, generalized $(\alpha, \beta)$-metrics, spherically symmetric Finsler metrics in $\mathbb{R}^n$, and other classes of Finsler metrics associated with the classification of Douglas spaces.

\section{Douglas curvature in some special classes of Finsler metrics}\label{DouglasSpecial}

The condition for a class of Finsler metrics to be of Douglas type has been a focal point in Finsler Geometry, prompting substantial research efforts in this field. Researchers have delved into identifying the precise conditions that distinguish Finsler metrics as Douglas metrics. This exploration has led to a deeper understanding of the intricate relationships between the class of Douglas metrics and different types of Finsler metrics. It aids in identifying gaps and new avenues for research within the realm of Douglas metrics.
We have divided section \ref{DouglasSpecial} into two following subsections, "Douglas Curvature and $(\alpha, \beta)$-metrics" and "Beyond $(\alpha, \beta)$-metrics".

\subsection{Douglas curvature and $(\alpha, \beta)$-metrics}
One of the most prominent classes of Finsler metrics are $(\alpha, \beta)$-metrics, which provide a connection to applications of Finsler metrics in various fields. Additionally, they are computationally tractable, allowing for the exploration and extension of their properties to broader classes of Finsler spaces. In this subsection, we begin by reviewing research on the topic of Douglas curvature, focusing specifically on spaces with $(\alpha, \beta)$-metrics.

In the early investigations focusing on metrics within the Douglas class, \cite{P*Finsler} revealed that every $P^*$-Randers metric with dimensions greater than three and vanishing Douglas curvature is a Riemannian metric. Additionally, a previous study \cite{P*FinslerIzumi} demonstrated that $C$-reducible $P^*$-Finsler metrics with dimensions exceeding three exhibit vanishing Douglas curvature. In the context of Finsler geometry, the study of $C$-reducible $P^*$-Finsler metrics in dimensions at most three remains an unexplored area within this paper. However, it is worth noting that a subsequent study, \cite{2-KropinaDouglas}, built upon the research by Matsumoto \cite{Matsomoto1998}, has established that any Kropina space of dimension two belongs to the class of Douglas metrics.

As a notable study on $(\alpha, \beta)$-metrics belonging to the class of Douglas metrics, the work \cite{ClassofDouglas} characterizes such metrics on an open subset $U \subset \mathbb{R}^n$ ($n \geq 3$) under the following conditions, with $b=\|\beta_x\|_{\alpha}$.
\ben
\item{} The 1-form $\beta$ is not parallel with respect to the Riemannian metric $\alpha$.
\item{} The $(\alpha, \beta)$-metric $F$ is not of Randers type.
\item{} $db \neq 0$ everywhere or $b = constant$ on $U$.
\een
A particularly interesting area of research is the study of Randers metrics that belong to the class of Douglas metrics. The paper \cite{DouglasRanders} has made remarkable contributions to this field by investigating Douglas-Randers manifolds with vanishing stretch curvature. Their research has shown that a Douglas-Randers manifold reduces to a Berwald manifold if and only if its stretch tensor vanishes. In this direction, the study \cite{Douglas-Squar} has explored Douglas-Square metrics with vanishing mean stretch curvature. Specifically, it has been shown that every Douglas-square metric ($n\geq 3$) is $R$-quadratic if and only if it is a Berwald metric.

The paper by \cite{SingularDouglasPF} investigates a specific class of $(\alpha, \beta)$-metrics that belong to the Douglas type, focusing on singular Finsler metrics. These metrics, which include Kropina metrics and $m$-Kropina metrics, have numerous applications in real-world scenarios. Notably, these metrics have been characterized by the conditions of being Douglasian and locally projectively flat in dimensions $n \geq 3$ and $n=2$, as outlined in \cite{SingularDouglasPF} and \cite{2-DimDouglasPF}, respectively. The main metrics studied in these two papers are of the form $\phi(s)=c s+s^m\varphi(s)$, where $c$ and $m$ are constant with $m \neq 0, 1$ and $\phi(s)$ is a $C^\infty$ function on a neighborhood of $s = 0$ satisfying $\phi(0) = 1$. The conditions under which these metrics belong to the class of Douglas metrics have been thoroughly investigated.

Other special classes of $(\alpha, \beta)$-metrics that belong to the Douglas type have been the focus of various research studies. For instance, the work by \cite{SingularSquareDouglas} delves into singular square metrics with vanishing Douglas curvature. Additionally, \cite{Kumar} has explored the conditions that determine when an $(\alpha, \beta)$-metric is classified as Douglas, applying this criterion to specific $(\alpha, \beta)$-metrics to investigate their properties when they fall under the category of Douglas metrics.

\subsection{Beyond $(\alpha, \beta)$-Metrics: Insights into Douglas Curvature}
The investigation of Douglas curvature in Finsler geometry is not limited to $(\alpha, \beta)$-metrics. Researchers have explored a wide range of Finsler metric classes, including generalized $(\alpha, \beta)$-metrics, as an extension of the existing studies focused on $(\alpha, \beta)$-metrics. The papers \cite{general(alphabeta)metricDouglas} and \cite{DouglasG(aB)} represent separate studies that have individually examined general $(\alpha, \beta)$-metrics of Douglas type. In \cite{general(alphabeta)metricDouglas}, H. Zhu identified a class of general $(\alpha, \beta)$-metrics with vanishing Douglas curvature under the condition that $\beta$ is closed and conformal with respect to $\alpha$. This means that the covariant derivatives of $\beta$ with respect to $\alpha$ satisfy $b_{i|j} = c a_{ij}$, where $c = c(x)\neq 0$ is a scalar function on $M$.
In \cite{DouglasG(aB)}, the authors remove the aforementioned condition and derive the equivalent equations for all Douglas general $(\alpha, \beta)$-metrics.
These works have introduced numerous novel Douglas metrics within this classification, contributing to the advancement of research in this area.

The paper \cite{SphDouglas}, following a similar approach to \cite{general(alphabeta)metricDouglas}, has identified the differential equation that defines spherically symmetric Finsler metrics with vanishing Douglas curvature. By solving this equation, the paper has successfully determined all the spherically symmetric Douglas metrics, providing a comprehensive understanding of this specific class of Finsler metrics. Additionally, the research includes numerous explicit examples to illustrate these findings, contributing significantly to the study of spherically symmetric Finsler metrics within the framework of Douglas curvature. In a specific scenario, the paper by \cite{SphericallyDouglas} has delineated spherically symmetric Douglas metrics with vanishing $S$-curvature and provided numerous explicit examples to illustrate these findings.

Beyond the well-established classes of Finsler metrics, the concept of vanishing Douglas curvature has been applied to various specialized classes of Finsler metrics. For instance, Ricci-flat Finsler metrics with the condition of vanishing Douglas curvature have been studied. Under the conditions of Ricci-flat Finsler metrics with vanishing Douglas curvature, a specific class of $(\alpha, \beta)$-metrics presented in \cite{ClassofDouglas} has been explored in two significant papers \cite{RicciFlatDpuglas} and \cite{Sevim}. Furthermore, \cite{RicciFlatDpuglas} has presented a significant theorem, demonstrating that an $(\alpha, \beta)$-metric of Randers type is a Ricci-flat Douglas metric if and only if it is a Berwald metric and $\alpha$ is Ricci-flat. This theorem provides a crucial characterization of Ricci-flat Douglas metrics within the context of Randers metrics.
\begin{thm}\cite{RicciFlatDpuglas}

Let $F$ be an $(\alpha, \beta)$-metric in the form $F = \sqrt{\alpha^2 + k\beta^3} + \varepsilon \beta$ on an $n$-dimensional manifold $M$ with $n\geq 2$, where $k$ and $\varepsilon \neq 0$ are constants. Then $F$ is a Ricci-flat Douglas metric if and only if $\beta$ is parallel with respect to $\alpha$ and the Riemann metric $\alpha$ is Ricci-flat. In this case, $F$ is a Berwald metric.
\end{thm}

The paper by \cite{HomogeneousDouglas(aB)} delves into an intriguing class of Finsler metrics known as homogeneous $(\alpha, \beta)$-metrics. It is demonstrated that

 \begin{thm} \cite{HomogeneousDouglas(aB)}

Let $F = \alpha \varphi(\frac{\beta}{\alpha})$ be a homogeneous $(\alpha, \beta)$-metric on $\frac{G}{H}$ . Then $F$ is a Douglas metric if and only if either $F$ is a Berwald metric or $F$ is a Douglas metric of Randers type.
\end{thm}

Another aspect of Finsler metrics explored in the context of vanishing Douglas curvature is the warped structures of Finsler metrics, as investigated in \cite{WarpedProductDouglas}. This study focuses on analyzing the warped structures of Finsler metrics to derive the differential equation that characterizes Finsler warped product metrics with vanishing Douglas curvature. Through the solution of this equation, the study identifies and describes all Finsler warped product Douglas metrics.

\subsection{\textbf{Douglas Curvature Review: Tables for Identifying Key Findings and Insights}}

In this subsection, we summarize the review of this section in the following table, providing a concise overview of the key findings and insights. This short summary enables researchers to rapidly pinpoint their areas of interest and navigate the research landscape surrounding Douglas curvature.
The following abbreviations were used in the tables bellow.
\vspace{8mm}
\begin{center}
\setlength\tabcolsep{0pt}
\begin{tabular}{|l |l| c|c|c|c|c|}
\hline
\textbf{Abbreviation} & \multicolumn{1}{|c|}{\textbf{Meaning}} \\
\hline
 S. & {Special}\\
 G. & {General}\\
 Lands. & Landsberg\\
 Sph. Symm. & Spherically Symmetric Finsler metrics in $\mathbb{R}^n$\\
 Tr. & Transformation\\
 Pr. & Projective\\
 Proj. & Projective\\
 Cur. & Curvature\\
 B & Berwald curvature\\
 L & Landsberg curvature\\
 E & E-curvature\\
 D & Douglas curcature\\
 S & S-curvature\\
 $\Sigma$ & Stretch curvature\\
 Douglas II & Douglas space of second kind\\
  L. Minkowski & Locally Minkowski\\
 Conf. Tr. & Conformal Transformatin\\
\hline
\end{tabular}
\end{center}

\newpage

\begin{center}
\setlength\tabcolsep{0pt}
 \begin{tabular}{|c|c|c|l|c|c|c|c|c|}
 \hline
 Res. Dir.& Key Findings & Ref &\multicolumn{1}{|c|}{Highlights}  \\
    \hline
\multirow{50}{*}{\begin{sideways}  Douglas curvature \end{sideways}}
 &  \multirow{10}{*}{G. Spaces} &\multirow{2}{*}{ \cite{Bacso}}& { \textcolor{green}{\textbf{-}} Specific conditions for certain classes of Finsler metrics to be of Douglas.}&  \\
 & & {} & {\textcolor{green}{\textbf{-}} $L=0$ + $D=0$  $\Leftrightarrow$ $B=0$.}\\
 \cline{3-4}

 & &  \multirow{3}{*}{\cite{Bacso2}} & {\textcolor{green}{\textbf{-}} Conditions for Special Douglas metrics}& \\
 && {} &{(Landsberg, G. $m$-Kropina ($\neq$ Kropina, Matsumoto) spaces) $\rightarrow$ Berwald.}& \\
 && {} &{\textcolor{green}{\textbf{-}} Douglas Surfaces with Cubic metric $\rightarrow$ L. Minkowski or Berwald.}  \\
 \cline{3-4}

 & &  \multirow{3}{*}{\cite{DShen}} & \\
 && {} &{\textcolor{green}{\textbf{-}} For non-Riemannian Douglas ($n \geq 3$):}& \\
 && {} &{ isotropic mean Berwald curvature $\Leftrightarrow$ relatively isotropic Landsberg curvature.} \\
 \cline{3-4}

 & & \multirow{2}{*}{ \cite{ReversibleDouglas}} & \\
 && {} &{\textcolor{green}{\textbf{-}} Compact Douglas + $\bar{E} = 0$ $\Rightarrow$ Berwald.}\\
 \cline{2-4}
 & & \multirow{2}{*}{\cite{P*FinslerIzumi}}&\\
 && {} &{\textcolor{green}{\textbf{-}} $C$-reducible + $P^*$-Finsler metrics ($n > 3$) $\Rightarrow$ $D=0$.}\\
 \cline{3-4}

 & & \multirow{2}{*}{\cite{P*Finsler}}&\\
 && {} &{\textcolor{green}{\textbf{-}} $P^*$-Randers metrics ($n > 3$) + $D=0$ $\Rightarrow$ Riemannian metric.}\\
 \cline{3-4}

 &  &  \multirow{2}{*}{ \cite{Kumar}} & \\
 && {} &{\textcolor{green}{\textbf{-}} Kropina metrics ($n = 2$) $\Rightarrow$ $D=0$ } \\
 \cline{3-4}

 &    &\multirow{5}{*}{ \cite{ClassofDouglas}} & \\
 && {} &{\textcolor{green}{\textbf{-}}  Characterizing $(\alpha, \beta)$-metrics on $U \subset \mathbb{R}^n$ ($n \geq 3$) to be Douglas when:}& \\
 && {} &{ $\diamond$ $\beta$ is not parallel w. r. t $\alpha$,} & \\
 && {} &{ $\diamond$ $F$ is not Randers,}& \\
 && {} &{ $\diamond$ $db \neq 0$ everywhere or $b = constant$ on $U$, with $b=\|\beta_x\|_{\alpha}$.}\\
 \cline{3-4}

 &  &  \multirow{2}{*}{ \cite{DouglasRanders}} &  \\
 && {}  & { \textcolor{green}{\textbf{-}} For Douglas-Randers:  $B=0$ $\Leftrightarrow$ $\Sigma=0$}\\
 \cline{3-4}
 &\multirow{1}{*}{$ (\alpha, \beta)- metrics $} &
 \multirow{2}{*}{\cite{Douglas-Squar}}&   \\
 && {}  & { \textcolor{green}{\textbf{-}} For Douglas-square metric ($n\geq 3$): $R$-quadratic $\Leftrightarrow$ $B=0$}\\
 \cline{3-4}

 &  &\multirow{3}{*}{\cite{SingularDouglasPF}} &  \\
 && {}  & { \textcolor{green}{\textbf{-}} ($n \geq3$) Investigating the conditions under which} & \\
 && {}  & { $F=\alpha \phi(s)$ with $\phi(s)=c s+s^m\varphi(s)$ belong to the class of Douglas metrics} \\
 \cline{3-4}

 &  &\multirow{2}{*}{\cite{2-DimDouglasPF}} &  \\
 && {}  & { \textcolor{green}{\textbf{-}} The significance of $n=2$ in \cite{SingularDouglasPF}}\\
 \cline{3-4}

 &  &\multirow{2}{*}{\cite{SingularSquareDouglas}} &  \\
 && {}  & { \textcolor{green}{\textbf{-}} Studying singular square metrics with vanishing Douglas curvature}\\
 \cline{3-4}

 &   &\multirow{3}{*}{\cite{Kumar}} &  \\
 && {}  & { \textcolor{green}{\textbf{-}} Investigating $(\alpha, \beta)$-metrics with Douglas curvature,}& \\
 && {} &{ applying this criterion to specific $(\alpha, \beta)$-metrics } \\
 \cline{2-4}
 &\multirow{4}{*}{$ G. (\alpha, \beta)- metric$}  &\multirow{2}{*}{\cite{general(alphabeta)metricDouglas}} & \\
 && {}  & { \textcolor{green}{\textbf{-}} Identifying general $(\alpha, \beta)$-metrics with $D=0$ under the condition  $b_{i|j} = c a_{ij}$. } \\
 \cline{3-4}

 & &\multirow{2}{*}{\cite{betaConformal}}&\\
 && {}  & { \textcolor{green}{\textbf{-}}  Identifying general $(\alpha, \beta)$-metrics with $D=0$ (without the condition $b_{i|j} = c a_{ij}$). }\\
 \cline{2-4}
 \parbox{0.4cm}&\multirow{4}{*}{ Sph. symm}  &\multirow{2}{*}{\cite{SphDouglas}} & \\
 && {}  & { \textcolor{green}{\textbf{-}} Identifying the diff. equation for  Sph. Symme. with $D=0$. } \\
 \cline{3-4}

 &  &\multirow{2}{*}{\cite{SphericallyDouglas}} & \\
 && {}  & { \textcolor{green}{\textbf{-}}  Identifying Sph. Symm. with $S=0$ and $D=0$.}\\
 \cline{2-4}
 \parbox{0.4cm}&\multirow{9}{*}{Other S. metrics} &\multirow{4}{*}{\cite{RicciFlatDpuglas}}& \\
 && {}  & { \textcolor{green}{\textbf{-}} For Randers metrics; $F$ is a Ricci-flat Douglas $\Leftrightarrow$ $B=0$ and $\alpha$ is Ricci-flat} & \\
 && {}  & { \textcolor{green}{\textbf{-}} $F = \sqrt{\alpha^2 + k\beta^3} + \varepsilon \beta$ ($n\geq 2$), (constants $k$ and $\varepsilon \neq 0$);}& \\
 && {}  & {$F$: Ricci-flat Douglas $\Leftrightarrow$ $\beta$ : parallel w. r. t. $\alpha$ and $\alpha$: Ricci-flat.}\\
 \cline{3-4}

 &  &\multirow{3}{*}{\cite{HomogeneousDouglas(aB)}} & \\
 && {}  & { \textcolor{green}{\textbf{-}} homogeneous $(\alpha, \beta)$-metric on $\frac{G}{H}$;  }& \\
 && {}  & {$F$: Douglas $\Leftrightarrow$ either $F$: Berwald or $F$: Randers Douglas metric.}\\
 \cline{3-4}

 &  &\multirow{2}{*}{\cite{WarpedProductDouglas}} & \\
 && {}  & { \textcolor{green}{\textbf{-}} Finding equation that characterizes Finsler warped product metrics + $D=0$.} \\
 \hline
\end{tabular}
 \end{center}

\section{Transformations in Douglas Spaces: A Comprehensive study}\label{Transformation}

This section delves into the world of metric transformations within Douglas spaces, reviewing the papers and researches that have contributed to a detailed exploration of the various changes and properties that characterize these spaces. By examining the existing literature, we will uncover the intricate relationships between transformations and the fundamental properties of Douglas metrics, allowing us to gain a deeper understanding of how these metrics behave under different changes in Finsler spaces. This analysis will also enable us to identify potential gaps that may exist in the current understanding of these transformations, providing a foundation for future research and exploration in this area.

\subsection{Conformal Transformations and Their Impact on Douglas Spaces}

This section provides a comprehensive overview of the significant contributions made to the study of conformal transformations and their effects on Douglas spaces. The studies have investigated various aspects, including the conditions under which vanishing Douglas curvature becomes conformally closed, the properties of conformal transformations of difference tensors, and the conformal invariant tensors associated with these difference tensors. Additionally, the section highlights the findings of specific papers that have examined the conformal transformations of specialized Finsler metrics of Douglas type, including the conditions under which they maintain their Douglas properties.

An early study that delved into transformations of Douglas metrics could be the work by \cite{Confclose}. This research focused on conformally closed Berwald and Douglas spaces, particularly investigating the conditions under which vanishing Douglas curvature becomes conformally closed within particular classes of Finsler metrics, such as Randers or Kropina metrics. Prior to that, the paper by \cite{ConfDiffTensors} explored the properties of conformal transformations of difference tensors in Finsler spaces with $(\alpha, \beta)$-metrics. Additionally, it introduced conformal invariant tensors in Finsler spaces with an $(\alpha, \beta)$-metric that are associated with these difference tensors.

Subsequent to the prior investigation, various research endeavors have been devoted to exploring the conformal transformations of specialized Finsler metrics of Douglas type, with a focus on determining the conditions under which they maintain their Douglas properties. In a contribution, the paper by \cite{ConfDouglas(aB)} examined the metric $F=c_1\alpha+c_2 \beta+\frac{\beta^2}{\alpha}$, where $c_2\neq 0$, and established that, under the assumption $\alpha^2 \not\equiv 0$ ($mod$ $\beta$), a Douglas space characterized by this metric undergoes a conformal transformation to another Douglas space if and only if the transformation is homothetic.

Similarly, in the related investigations, \cite{ConfDouglasAppeox}, \cite{ConfDouglasSpecial(aB)} and \cite{Conf(aB)} examined the second approximate Matsumoto metric given by $F=\alpha+\beta+\frac{\beta^2}{\alpha}+\frac{\beta^3}{\alpha^2}$, $F=\alpha-\frac{\beta^2}{\alpha}$ and $F=\kappa(\alpha+\beta)+\varepsilon\frac{\beta^2}{\alpha}$, (with non-zero constants $\kappa$ and $\varepsilon$), respectively. These studies aimed to identify the conditions under which these metrics belong to the Douglas type class. Furthermore, the authors determined the criteria for a conformally transformed Douglas space with the specified metrics to remain a Douglas space.

The paper \cite{ConfDouglas} presents a remarkable result, showing that two-dimensional conformally related Douglas metrics are Randers metrics. This finding highlights the special nature of Douglas metrics in dimension two.

It is important to note that Douglas surfaces and Finsler surfaces in general can exhibit different behaviors in various problems. For instance, in the case of Douglas surfaces, the work of \cite{DouglasSurface} is noteworthy.
In \cite{DShen}, Chen and Shen proved that on a Douglas manifold with dimension $n \geq 3$, a Finsler metric has isotropic mean Berwald curvature if and only if it has relatively isotropic Landsberg curvature. However, this fact remained unsolved for the case of two-dimensional Finsler manifolds until the following theorem was proved.

\begin{thm}

Let $(M, F)$ be a Douglas surface. Then the following are equivalent.
\ben
\item{} $F$ has isotropic mean Berwald curvature.
\item{} $F$ has relatively isotropic Landsberg curvature.
\een
\end{thm}

General $(\alpha, \beta)$-metrics have been studied under conformal changes with the condition of Douglas type. \cite{ConfG(aB)} establishes the same result as the previous papers, which states that if a general $(\alpha, \beta)$-metric of the form $F=\alpha \varphi(b^2, s)$ is a Douglas metric, then its conformal transform $\bar{F} = e^{\sigma}F$ is also a Douglas metric if and only if the conformal transformation is a homothety, under suitable conditions on the conformal factor $\varphi$.
Specifically, \cite{ConfG(aB)} proves the following theorem:

\begin{thm}

If the general $(\alpha, \beta)$-metric $F=\alpha \varphi(b^2, s)$ is a Douglas metric, then its conformal transform $\bar{F} = e^{\sigma}F$ is also a Douglas metric if and only if the conformal transformation is a homothety, under suitable following conditions on the function $\varphi$.
\[
\varphi_1 \varphi_2-\varphi \varphi_{12} \neq 0,
\]
where $\varphi_1=\frac{\partial \varphi}{\partial b^2}$, $\varphi_2=\frac{\partial \varphi}{\partial s}$, $\varphi_{12}=\frac{\partial^ \varphi}{\partial b^2 \partial s}$ and $\sigma=\sigma(b^2)$.
\end{thm}

\subsubsection{\textbf{Conformal Transformations and Douglas spaces of the second kind}}
However, in this study, the focus was on the condition of Douglas of the second kind rather than the Douglas type I. Y. Lee \cite{DouglasSecondType} recently explored Douglas spaces of the second kind and established criteria for a Finsler space with Matsumoto metric to be a second kind Douglas space. In \cite{ConfDouglasSecondG(aB)}, it was demonstrated that a second kind Douglas space with a generalized form of $(\alpha, \beta)$-metric $F$ can be conformally transformed into another second kind Douglas space.

There is a line of researches that has focused on studying the conformal changes of different classes of Finsler metrics that are Douglas spaces of the second kind. These studies aim to determine the conditions under which a Douglas space of the second kind with a specific class of Finsler metric is conformally transformed into another Douglas space of the second kind.

In \cite{ConfSecondDouglas(ab)2015}, the authors furthermore found the conditions under which the conformal change of a Finsler space with Matsumoto and generalized Kropina metrics results in a Douglas space of the second kind. Building upon this, \cite{ConfSecondDouglasSpecial(aB)} continued the investigation and determined the conditions for a Finsler space with the special $(\alpha,\beta)$-metric $F=\beta+\frac{\beta^2}{\alpha}$ to be conformally transformed into a Douglas space of the second kind.

In a similar context, both \cite{ConftransDouglasSecond(aB)} and \cite{Confchange(aB)} have illustrated that Douglas spaces of the second kind possess specific metrics that undergo conformal transformations, resulting in Douglas spaces of the second kind.
Specifically, \cite{ConftransDouglasSecond(aB)} examined the metric $F=\kappa(\alpha+\beta)+\varepsilon\frac{\beta^2}{\alpha}$, where $\kappa$ and $\varepsilon$ are non-zero constants, while \cite{Confchange(aB)} investigated the metric $F=c_1\alpha+c_2\beta+c_3\frac{\alpha^2}{\beta}+\frac{\beta^2}{\alpha}$ on a manifold $M$ with dimension $n>2$, where $c_i$, ($i=1, 2, 3$), are real constants. Both studies concluded that a Douglas space of the second kind with these metrics undergoes conformal transformations to become another Douglas space of the second kind.

\subsection{Exploring Beyond Conformal Transformations}

This section presents the investigations that have done on the various types of transformations and their effects on the properties and structures of Douglas spaces, extending the analysis beyond the realm of conformal transformations to uncover new insights and advancements in the field.

An early paper, such as \cite{RandersChange(aB)}, explores the concept of a special Randers change involving $(\alpha, \beta)$-metrics of Douglas type, demonstrating that they retain the Douglas type property under this transformation, and conversely.

The notion of $C$-conformal transformation, a specific type of conformal transformation, was introduced by Hashiguchi \cite{Hashiguchi}. In \cite{C-Conformal}, the authors consider the infinitesimal $C$-conformal motion of certain Finsler spaces, including Douglas spaces.

$\beta$-change, generalized $\beta$-change, (Almost) $\beta$-change and (Almost) $\beta$-conformal change are the changes that have been considered in \cite{betaChangeSpecial}, \cite{betaConformal} and \cite{AlmostbetaChange}, respectively. The paper \cite{betaConformal} explores the projective generalized $\beta$-change, represented by the form $\bar{F}=f(e^{\sigma}F, \beta)$, applied to a Minkowskian Finsler metric $F$, resulting in the vanishing of the Douglas curvature. In papers \cite{betaChangeSpecial} and \cite{AlmostbetaChange}, the conditions for Douglas $(\alpha, \beta)$-metrics $F$ and spherically symmetric Finsler metrics in $\mathbb{R}^n$, denoted as $F=u\varphi(r,s)$, to remain invariant under the $\beta$-change, (Almost) $\beta$-change and (Almost) $\beta$-conformal change are investigated. These transformations are defined as $\bar{F}=F+\beta$, $\bar{\varphi}=\varphi+h(s)$ and $\bar{\varphi}=e^{\sigma}\varphi+h(s)$ for the (Almost) $\beta$-change and (Almost) $\beta$-conformal change, respectively. The notations $s$, $u$, $v$, $\varphi$, and other terms pertaining to spherically symmetric Finsler metrics are defined as in \cite{betaChangeSpecial} and \cite{AlmostbetaChange}.

Various transformations of Finsler metrics are explored in the papers \cite{FirstApproxExpo}, \cite{MatsumotoChangeh-vector}, and \cite{ConfKropRanders}. These studies investigate the conditions under which specific Finsler metrics of Douglas type maintain their Douglas properties through transformations such as the First Approximate Exponential change of Finsler Metric $F$ given by $\bar{F}= F+ \beta+\frac{\beta^2}{2F}$, the Matsumoto change of Finsler Metric $F$ using the $h$-vector represented as $\bar{F}=\frac{F^2}{F-\beta}$ with $h$-vector $b_i$, and the conformal Kropina-Randers change of $(\alpha, \beta)$-metric $F$ denoted as $\bar{F}= e^{\sigma}(\frac{F^2}{\beta}+\beta)$.
It is noteworthy that the concept of the $h$-vector was introduced by Izumi in \cite{Izumi} while investigating the conformal transformation of Finsler spaces. The $h$-vector $b$ with components $b_i$ is defined as $v$-covariantly constant with respect to Cartan's connection, satisfying $FC{^r}_{jk}b_r=\rho h_{jk}$, where $\rho$ is a scalar function of $x$, $y$, and $h_{jk}$ represents the components of the angular metric tensor.

\subsection{\textbf{Transformations in Douglas Spaces Review: Tables for Identifying Key Findings and Insights}}

In this subsection, we present a summary of this section in the following table, offering a concise overview of the key findings and insights. This brief summary allows researchers to quickly identify their areas of interest and navigate the research landscape related to transformations in Douglas spaces.
\begin{center}
\setlength\tabcolsep{0pt}
 \begin{tabular}{|c|c|c|l|c|c|c|c|c|}
 \hline
Res. Dir.& Key Findings & Ref &\multicolumn{1}{|c|}{Highlights}  \\
    \hline
\multirow{35}{*}{\begin{sideways}  Transformations and Douglas cur. \end{sideways}}

 & & \multirow{2}{*}{\cite{ConfDiffTensors}} & {\textcolor{green}{\textbf{-}} The conformal transformations of difference tensors }& \\
 && {} &{in Finsler spaces with $(\alpha, \beta)$-metrics.}  \\

 \cline{3-4}
 & & \multirow{2}{*}{\cite{Confclose}}& { \textcolor{green}{\textbf{-}} Investigating the conditions for conformally closed } &  \\
 & & {} & {Douglas metrics within particular Finsler metrics.}\\
 \cline{3-4}

 & &\multirow{3}{*}{\cite{ConfDouglas(aB)}} & \\
 && {} &{\textcolor{green}{\textbf{-}} For $F=c_1\alpha+c_2 \beta+\frac{\beta^2}{\alpha}$, with $c_2\neq 0$, $\alpha^2 \not\equiv 0$ ($mod$ $\beta$):}& \\
 && {} &{ Douglas metric is conformally closed $\Leftrightarrow$ it is homothetic.} \\
 \cline{3-4}

 & \multirow{6}{*}{ Conf. Tr. and Douglas $(\alpha, \beta)$-metrics} &  \multirow{3}{*}{\cite{ConfDouglasAppeox}} & \\
 && {} &{\textcolor{green}{\textbf{-}} For $F=\alpha+\beta+\frac{\beta^2}{\alpha}+\frac{\beta^3}{\alpha^2}$:}& \\
 && {} &{ Douglas metric is conformally closed $\Leftrightarrow$ it is homothetic.} \\
 \cline{3-4}

 & & \multirow{3}{*}{ \cite{ConfDouglasSpecial(aB)}} & \\
 && {} &{\textcolor{green}{\textbf{-}} For $F=\alpha-\frac{\beta^2}{\alpha}$:}& \\
 && {} &{ Douglas metric is conformally closed $\Leftrightarrow$ it is homothetic.} \\
 \cline{3-4}

 & &  \multirow{3}{*}{\cite{Conf(aB)}} & \\
 && {} &{\textcolor{green}{\textbf{-}} For $F=\kappa(\alpha+\beta)+\varepsilon\frac{\beta^2}{\alpha}$, (non-zero constants $\kappa, \varepsilon$):}& \\
 && {} &{ Douglas metric is conformally closed $\Leftrightarrow$ it is homothetic.} \\
 \cline{3-4}

 & &\multirow{2}{*}{\cite{ConfDouglas}} & \\
 && {} &{\textcolor{green}{\textbf{-}} 2-dimensional conformally related Douglas $\Rightarrow$ Randers.}\\
 \cline{2-4}
 & \multirow{4}{*}{ Conf. Tr. and G. Douglas $(\alpha, \beta)$-metrics}
 & \multirow{4}{*}{\cite{ConfG(aB)}} & \\
 && {} &{\textcolor{green}{\textbf{-}}  For G. $(\alpha, \beta)$-metrics $F=\alpha \varphi(b^2, s)$ }& \\
 && {} &{ Douglas metric is conformally closed $\Leftrightarrow$ it is homothetic }& \\
 && {} &{ under condition $\varphi_1 \varphi_2-\varphi \varphi_{12} \neq 0$.} \\
 \cline{2-4}
 & &\multirow{4}{*}{\cite{DouglasSecondType}}&\\
 && {} &{\textcolor{green}{\textbf{-}} Finding conditions for Matsumoto metrics to be}&\\
 && {} &{ a Douglas II spaces.}\\
 \cline{3-4}

 & \multirow{6}{*}{ Conf. Tr. and Douglas II } &\multirow{3}{*}{ \cite{ConfDouglasSecondG(aB)}}&\\
 && {} &{\textcolor{green}{\textbf{-}}  Conf. invariant Douglas II with G. $(\alpha, \beta)$-metric. }& \\
 && {} &{\textcolor{green}{\textbf{-}} Conf. invariant Douglas II with Randers and }& \\
 && {} &{G. $m$-Kropina metrics.}\\
 \cline{3-4}

 & &\multirow{3}{*}{\cite{ConfSecondDouglas(ab)2015}} & \\
 && {}  & { \textcolor{green}{\textbf{-}} Identifying conditions under which a Douglas II is}& \\
 && {} &{Conf. closed with Matsumoto and G. Kropina metrics. } \\
 \cline{3-4}
 \hline
\end{tabular}
\end{center}
\begin{center}
\setlength\tabcolsep{0pt}
 \begin{tabular}{|c|c|c|l|c|c|c|c|c|}
 \hline
Res. Dir. & Key Findings & Ref &\multicolumn{1}{|c|}{Highlights } \\
    \hline
 & &\multirow{3}{*}{\cite{ConfSecondDouglasSpecial(aB)}} & \\
 && {}  & { \textcolor{green}{\textbf{-}} Identifying conditions under which a Douglas II is}& \\
 && {} &{Conf. closed with $F=\beta+\frac{\beta^2}{\alpha}$. } \\
 \cline{3-4}

 &\multirow{6}{*}{ Conf. Tr. and Douglas II }  &\multirow{4}{*}{\cite{ConftransDouglasSecond(aB)}} & \\
 && {}  & { \textcolor{green}{\textbf{-}} Identifying conditions under which a Douglas II is}& \\
 && {} &{ Conf. closed with $F=\kappa(\alpha+\beta)+\varepsilon\frac{\beta^2}{\alpha}$. } \\
 \cline{3-4}

 & &\multirow{3}{*}{\cite{Confchange(aB)}} & \\
 && {}  & { \textcolor{green}{\textbf{-}} Identifying conditions under which a Douglas II ($n>2$) is}& \\
 && {} &{Conf. closed with $F=c_1\alpha+c_2\beta+c_3\frac{\alpha^2}{\beta}+\frac{\beta^2}{\alpha}$. } \\
 \cline{2-4}

 && \multirow{3}{*}{\cite{RandersChange(aB)} }& \\
 && {}  & { \textcolor{green}{\textbf{-}} Finding conditions for Finsler spaces changed by }& \\
 && {} &{ a special Randers change to be of Douglas type. } \\
 \cline{3-4}

\multirow{3}{*}{\begin{sideways}  Transformations and Douglas cur. \end{sideways}}

 &&\multirow{3}{*}{\cite{C-Conformal}} & \\
 && {}  & { \textcolor{green}{\textbf{-}}  Studying $C$-conformal change of certain Finsler spaces,}& \\
 && {} &{including Douglas spaces.}\\
 \cline{3-4}
 & & \multirow{3}{*}{\cite{betaConformal}} & \\
 && {}  & { \textcolor{green}{\textbf{-}} Proj. G. $\beta$-change of Minkowski metric $F$, $\bar{F}=f(e^{\sigma}F, \beta)$, } & \\
 && {} &{ $\Rightarrow$ $D=0$ } \\
 \cline{3-4}

 & Beyond Conf. Tr. and Douglas Spaces
 &\multirow{3}{*}{\cite{betaChangeSpecial}} & \\
 && {}  & { \textcolor{green}{\textbf{-}}  Investigating invariant Douglas spaces with $(\alpha, \beta)$-metrics }& \\
 && {} &{under the $\beta$-change.}\\
 \cline{3-4}

 & &\multirow{3}{*}{\cite{AlmostbetaChange}} & \\
 && {}  & { \textcolor{green}{\textbf{-}}  Investigating invariant Douglas spaces with Sph. Symm. }& \\
 && {} &{under the (Almost) $\beta$ and (Almost) $\beta$-conformal changes.}\\
 \cline{3-4}

 & & \multirow{2}{*}{ \cite{ConfDiffTensors}} & {\textcolor{green}{\textbf{-}} Finding conditions for invariant Douglas spaces with }& \\
 && {} &{the change $\bar{F}= F+ \beta+\frac{\beta^2}{2F}$. } \\
 \cline{3-4}

 & & \multirow{3}{*}{\cite{MatsumotoChangeh-vector}} & {\textcolor{green}{\textbf{-}} Finding conditions for invariant Douglas spaces with}& \\
 && {} &{Matsumoto change of $F$ using the $h$-vector $b_i$, $\bar{F}=\frac{F^2}{F-\beta}$. } \\
 \cline{3-4}

 & & \multirow{2}{*}{\cite{ConfKropRanders}} & {\textcolor{green}{\textbf{-}} Finding conditions for invariant Douglas spaces with}& \\
 && {} &{Conf. Kropina-Randers Tr. of $(\alpha, \beta)$-metric, $\bar{F}= e^{\sigma}(\frac{F^2}{\beta}+\beta)$. } \\
 \hline
\end{tabular}
\end{center}
\newpage
\section{The generalization of Douglas curvature: Expanding the Horizons of Douglas Spaces}\label{Generalization}

After the seminal work by Sakaguchi \cite{Sakaguchi} delving into the well-known general Douglas space, now recognized as General Douglas-Weyl ($GDW$)-metrics, B\'{a}cs\'{o} demonstrated in \cite{BacsoPapp} that this category remains invariant under projective changes. Furthermore, the research in \cite{Sakaguchi} established that all Weyl metrics are encompassed within this broader class of general Douglas spaces. Subsequently, an increasing number of studies have focused on this intriguing class of Finsler metrics.

In \cite{subclassGDW}, a subclass of $GDW$-metrics was introduced, characterized by a constant Douglas curvature along Finslerian geodesics, revealing that this subclass includes $R$-quadratic Finsler metrics. Additionally, the paper investigates the conditions under which this specific class of Finsler metrics transitions into the category of Landsberg metrics.

In the research presented in \cite{specialclassGDW}, the authors further investigate the subclass of $GDW$-metrics characterized by the property $D_j{^i}_{kl|m}y^m=0$. They demonstrate that for each Landsberg metric within this subset of Finsler metrics, the condition $\bar{E}_{jkl} =E_{jk|l}= 0$ holds if and only if $H = 0$. Additionally, they establish that any Finsler metric exhibiting non-zero isotropic flag curvature in this specific subclass of Finsler metrics is Riemannian if and only if $\bar{E} = 0$.

A unique category of Finsler metrics, known as $D$-recurrent Kropina space, is examined to be part of this generalized Douglas space. The investigation into this specific subclass of $GDW$-metrics is detailed in \cite{Ceyhan1}. Moreover, it is proved that

\begin{thm}\cite{Ceyhan1}

For a $D$-recurrent Kropina space $F^n$ ($n>2$) with weak Berwald metric, the followings are equivalent,
\ben
\item{} $F$ is a Berwald space,
\item{} $F$ is a generalized Douglas space.
\een
\end{thm}

The investigation into $GDW$-metrics was extended in \cite{Generalized Douglas-Weyl Spaces}, revealing that the category of $R$-quadratic Finsler spaces is a proper subset of the broader class of generalized Douglas-Weyl spaces. Furthermore, it was demonstrated that all generalized Douglas-Weyl spaces exhibiting zero Landsberg curvature also exhibit zero non-Riemannian quantity $H$, expanding upon a previously established result limited to $R$-quadratic metrics. The paper asserts that this advancement also results in an extension of the well-known Numata's Theorem, which states that every Landsberg metric $F$ ($n \geq 3$) with non-zero scalar flag curvature is a Riemannian metric with constant sectional curvature \cite{Numata}.

The paper \cite{Emamian} explores a specific class of $(\alpha, \beta)$-metrics with vanishing $S$-curvature and demonstrates that they belong to the category of $GDW$-metrics. In particular,

\begin{thm}\cite{Emamian}

Let $F=\alpha \varphi (s)$, $s=\frac{\beta}{\alpha}$, be an $(\alpha, \beta)$-metric on a manifold $M$ of dimension $n\geq 3$. Suppose that
\[
F \neq c_s(\frac{\beta}{\alpha})^{\frac{c_2}{1+c_2}}\big(c_1\frac{\beta}{\alpha}+c_2+1\big)^{\frac{1}{1+c_2}}, \quad and \quad
F\neq d_1\sqrt{\alpha^2+d_2 \beta^2}+ d_3 \beta.
\]
where $c_1$, $c_2$, $c_3$, $d_1$, $d_2$ and $d_3$ are real constants. Let $F$ has vanishing $S$-curvature. Then $F$ is a $GDW$-metric if and only if it is a Berwald metric.
\end{thm}

To further investigate the unique category of $GDW$-metrics, \cite{GDWRanders} delineates $GDW$-Randers metrics by their Zermelo navigation data. Additionally, the paper introduces a set of $GDW$-Randers metrics that deviate from being $R$-quadratic. However, in \cite{PFRanders}, the following theorem has been established.

\begin{thm}\cite{PFRanders}

Let $F = \alpha + \beta$ be a Randers metric on an $n$-dimensional manifold $M$. Then $F$ is a $GDW$-metric if and only if
\[
s_{ij|k}=\frac{1}{n-1}(a_{ik}s{^m}_{j|m}-a_{jk}s{^m}_{i|k}),
\]
where "$\mid$" denotes the covariant derivative with respect to $\alpha$.
\end{thm}

Note that $s_{ij}=\frac{1}{2}(b_{i|j}-b_{j|i})$.

The intriguing paper \cite{GDWLie} investigates a Lie group $G$ equipped with a left-invariant Randers metric $F$ and establishes a connection between the projective geometry of $(T G, F^c)$ and $(T G, F^v)$ and the projective geometry of $(G, F)$, where $F^c$ and $F^v$ represent the complete and vertical lift of $F$, respectively. In particular,

\begin{thm}\cite{GDWLie}

Let $G$ be an $n$-dimensional Lie group equipped with a left-invariant Randers metric $F(x, y) =\sqrt{\alpha_x(y, y)} + \alpha_x(U, y)$ defined by the underlying left-invariant Riemannian metric $\alpha$ and the left-invariant vector field $U$ such that $\|U\|_{\alpha} < 1$. Then the followings hold
\ben
\item[(i)]{} $F^c$ is a $GDW$-metric on $TG$  if and only if the following conditions are satisfied.
\[
<U, [\nabla{}_{X} Y-\frac{1}{2}[X,Y], Z]>
\]
\[
=\frac{1}{2n-1}<X,Y>\sum \big<U, 2[\nabla_{X_i}X_i, Z]+[X_i, \nabla_{X_i}Z] \big>,
\]
\[
<U, \big[[X,Z], Y\big]>+<U,\big[[X,Y],Z\big]>=0,
\]
where $X$, $Y$ and $Z$ are arbitrary left-invariant vector fields, and $\{X_i\}_{i=1}^n$ is an orthogonal basis of the Lie algebra $g$ with respect to $\alpha$
\item[(ii)]{} $F^v$ is a $GDW$-metric on $TG$  if and only if the following conditions are satisfied.
\[
\mathcal{C}(X,Y,Z)-\frac{1}{2}<U,[[X, Y], Z] + [Y, [X, Z]]>
\]
\[ = \frac{1}{2n-1} \{<X,Y><U, B(Z)>- <X,Z><U, B(Y)> \},
\]
\[
\mathcal{C}(X,Y,Z)+\frac{1}{2}<U,[Y, ad^*_ZX]>=\frac{1}{2n-1} <X,Y><U, B(Z)>,
\]
\[
\mathcal{C}(X,Y,Z)+\frac{1}{2}<U,[ad^*_XY,Z]>+ [Y, ad^*_XZ]=0,
\]
where $X, Y, Z \in \mathfrak{g}$.
\een
\end{thm}

All the symbols and terms used in the aforementioned theorem are defined and clarified in detail in the publication referenced as \cite{GDWLie}. Certainly, there exist additional extensions of Douglas metrics beyond $GDW$-metrics that warrant attention. In \cite{ReversibleDouglas}, a novel notion called reversible Douglas curvature is introduced, defining it as a Finsler metric where $D_j{^i}_{kl}(x,y)=D_j{^i}_{kl}(x,-y)$. The study reveals that a Finsler metric qualifies as a Douglas metric precisely when it exhibits reversible geodesics and reversible Douglas curvature. It has been demonstrated that Randers metrics must exhibit reversible Douglas curvature, while Kropina metrics are required to possess reversible geodesics, leading to the natural presence of reversible Douglas curvature in the latter. The concept of reversibility is detailed in \cite{ReversibleDouglas}, where a reflection map $\rho : (x,y) \in TM \rightarrow (x, -y) \in TM$ is used to define the reverse of the geodesic spray $G$ as $\bar{G}= -\rho^* G$. See also \cite{ReversibleDouglas} for geometrical and dynamical characterizations of geodesically reversible Finsler metrics. Consequently, the focus shifts towards the broader category of Finsler metrics featuring reversible Douglas curvature, a topic deserving significant attention. This prompts further exploration, as seen in the investigation of alternative metrics like square metrics, aiming to address the following two key inquiries.

- How extensive is the set of Finsler metrics with reversible Douglas curvature?

- What are the local construction of the Finsler metrics with reversible Douglas curvature expect for Randers metrics and Douglas metrics?

\begin{thm}\cite{ReversibleDouglas}

Let $F=\frac{(\alpha+\beta)^2}{\alpha}$  be a square metric on manifold $M$ of dimension $n>2$. Then $F$ has reversible Douglas curvature if and only if it is a Douglas metric.
\end{thm}

Another generalization of Douglas metrics is explored in \cite{GDouglas}, which investigates a novel class of Finsler metrics encompassing Douglas metrics as a special case, referred to as generalized Douglas metrics. The paper demonstrates that a generalized Douglas Finsler metric reduces to a Douglas metric if it is a $\tilde{B}$-stretch metric. A $\tilde{B}$-stretch metric is defined as a Finsler metric satisfying the condition
\[
\tilde{B}_j{^i}_{kl|m}=\tilde{B}_j{^i}_{km|l},
\]
where $\tilde{B}_j{^i}_{kl} ={B}_j{^i}_{kl|0}$ as defined in \cite{GDouglas}.
A generalization of Douglas metric which is some-how related to the previous paper, is presented in \cite{Sad1}. The paper constructs the new sub-classes of generalized Douglas-Weyl metrics and presents illustrative examples.  It studies a new quantity in Finsler geometry, so-called generalized Berwald projective Weyl ($GB\widetilde{W}$) curvature, which is a $C$-projective invariant and is the proper subset of the class of $GDW$-metrics.

In the realm of generalizing Douglas metrics, several studies have delved into the concept of other versions beyond the familiar natural generalized Douglas metrics. Notable examples include the works referenced in \cite{WD1}, \cite{WD2}, and \cite{WD3}. These papers introduce an intriguing projective invariant termed weakly Douglas Finsler metric, defined by $D_j{^i}_{kl}=T_{jkl}y^i$, where $T_{jkl}$ represents specific tensors. \\
In \cite{WD1}, it is demonstrated that every Randers manifold of dimension $n\geq 3$ is a weakly Douglas metric if and only if it is a Douglas metric. Additionally, it is established that every Kropina manifold is a weakly Douglas metric if and only if it is a Douglas metric, highlighting that every Kropina surface is a Douglas surface.

Two separate papers, \cite{WD2} and \cite{WD3}, individually prove that every weakly Douglas warped product Finsler metric is a Douglas metric. Furthermore, in \cite{WD2}, it is shown that every weakly Douglas warped product Finsler metric is Berwald if and only if it is Landsberg, while in \cite{WD3}, under specific conditions, a class of Finsler warped product metrics is locally projectively flat if and only if it exhibits scalar flag curvature.

\subsection{\textbf{The generalization of Douglas curvature: Tables for Identifying Key Findings and Insights}}

In this subsection, we offer a consolidated overview of the review in the table below, which emphasizes the key findings and insights. This concise summary enables researchers to swiftly pinpoint their areas of interest and navigate the research landscape concerning the generalization of Douglas curvature.
\newpage
\begin{center}
\setlength\tabcolsep{0pt}
 \begin{tabular}{|c|c|c|l|c|c|c|c|c|}
 \hline
Res. Dir. & Key Findings & Ref &\multicolumn{1}{|c|}{Highlights } \\
    \hline
\multirow{45}{*}{\begin{sideways}  The generalization of Douglas curvature \end{sideways}}
 &  &\multirow{3}{*}{\cite{Sakaguchi}} & \\
 && {}  & { \textcolor{green}{\textbf{-}}  Studying well-known G. Douglas space, $GDW$-metrics.}& \\
 && {} &{\textcolor{green}{\textbf{-}} Weyl metrics $\subseteq$ $GDW$-metrics.}\\
 \cline{3-4}

 & &\multirow{2}{*}{\cite{BacsoPapp}} & \\
 && {}  & { \textcolor{green}{\textbf{-}}  Invariance $GDW$-metrics under Proj. Tr.}\\
 \cline{3-4}

 & $GDW$-metrics
 &  \multirow{2}{*}{\cite{subclassGDW}} & {\textcolor{green}{\textbf{-}} Studying subclass of $GDW$-metrics + $D_j{^i}_{kl|0}=0$. }& \\
 && {} &{\textcolor{green}{\textbf{-}} Compact Finsler space + $D_j{^i}_{kl|0}=0$ and $\bar{E}= 0$ $\Rightarrow$ Lands.}\\
 \cline{3-4}

 & & \multirow{2}{*}{\cite{specialclassGDW}} & {\textcolor{green}{\textbf{-}} Finsler metric of isotropic flag Cur. ($\lambda\neq 0$) with $D_j{^i}_{kl|0}=0$: }& \\
 && {} &{Riemannian $\Leftrightarrow$ $\bar{E} = 0$.} \\
 \cline{3-4}

 & & \multirow{2}{*}{\cite{ConfKropRanders}} & {\textcolor{green}{\textbf{-}} Finding conditions for invariance Douglas spaces with}& \\
 && {} &{Conf. Kropina-Randers Tr. of $(\alpha, \beta)$-metric, $\bar{F}= e^{\sigma}(\frac{F^2}{\beta}+\beta)$. } \\
 \cline{2-4}
 & \multirow{6}{*}{other generalizations of Douglas Cur.} & \multirow{2}{*}{\cite{D-recurrent}} & \\
 && {}  & { \textcolor{green}{\textbf{-}}  Introducing $D$-recurrent Kropina space.}\\
 \cline{3-4}

 & \multirow{6}{*}{($D$-recurrent metrics)} &\multirow{4}{*}{\cite{Ceyhan1}} & \\
 && {}  & { \textcolor{green}{\textbf{-}}  For $D$-recurrent Kropina space($n>2$) + $E=0$: }& \\
 && {} &{$F$ is a Berwald space,}& \\
 && {} &{$F$ is a generalized Douglas-Weyl space.}\\
 \cline{2-4}

 & &\multirow{3}{*}{\cite{Generalized Douglas-Weyl Spaces}} & \\
 && {}  & { \textcolor{green}{\textbf{-}}  $GDW$-spaces + $L=0$ $\Rightarrow$ $H=0$. }& \\
 && {} &{($\Rightarrow$ An extension Numata Theorem).}\\
 \cline{3-4}

 &&\multirow{3}{*}{ \cite{Emamian}} & {\textcolor{green}{\textbf{-}} For some class of $(\alpha, \beta)$-metrics + $S=0$.}& \\
 && {} &{$F$ is a $GDW$-metric $\Leftrightarrow$ $B=0$. } \\
 \cline{3-4}

 & & \multirow{2}{*}{\cite{GDWRanders}} & {\textcolor{green}{\textbf{-}} Investigating the class of $GDW$-metrics}& \\
 && {} &{Introducing $GDW$-Randers $\neq$ $R$-quadratic.} \\
 \cline{3-4}

 & & \multirow{2}{*}{\cite{PFRanders}}& {\textcolor{green}{\textbf{-}} Randers metric is a $GDW$-metric if and only if: }& \\
 && {} & { $s_{ij|k}=\frac{1}{n-1}(a_{ik}s{^m}_{j|m}-a_{jk}s{^m}_{i|k})$. } \\
 & & \multirow{2}{*}{\cite{Generalized Douglas-Weyl Spaces}} & {\textcolor{green}{\textbf{-}} For some class of $(\alpha, \beta)$-metrics + $S=0$.}& \\
 && {} &{($\Rightarrow$ An extension Numata Theorem).}\\
 \cline{3-4}

 & & \multirow{2}{*}{\cite{PFRanders}} & {\textcolor{green}{\textbf{-}} Randers metric is a $GDW$-metric if and only if: }& \\
 && {} & { $s_{ij|k}=\frac{1}{n-1}(a_{ik}s{^m}_{j|m}-a_{jk}s{^m}_{i|k})$. } \\
 \cline{3-4}

 & \multirow{6}{*}{$GDW$ and S. Class of Finsler metrics}
 &\multirow{2}{*}{ \cite{GDWLie} }& {\textcolor{green}{\textbf{-}} Left-invariant Randers metric on a Lie group is }&\\
 && {} & {$GDW$-metric if and only if it satisfies specific conditions.} \\
 \cline{3-4}

 && \multirow{6}{*}{\cite{ReversibleDouglas}} & {\textcolor{green}{\textbf{-}} Introducing reversible Douglas Cur.}&\\
 && {} &  {\textcolor{green}{\textbf{-}} $F$ is a Douglas metric $\Leftrightarrow$ }& \\
 && {} & {it exhibits reversible geodesics and reversible Douglas Cur.}&\\
 && {} &  {\textcolor{green}{\textbf{-}} Investigating metrics with reversible Douglas Cur.}&\\
 && {} &  {to determine their extent and local construction.}&\\
 && {} &  {\textcolor{green}{\textbf{-}} a form square metric has reversible Douglas Cur. $\Leftrightarrow$ $D=0$. }\\
 \cline{3-4}

 && \multirow{3}{*}{\cite{Sad1}} & {\textcolor{green}{\textbf{-}} Presenting a new quantity, ($GB\widetilde{W}$) curvature,}&\\
 & & {} & { which is a $C$-projective invariant.}&\\
 && {} &  {\textcolor{green}{\textbf{-}} $GB\widetilde{W}(M)$ $\subseteq$ $GDW(M)$.  }\\
 \cline{3-4}

 && \multirow{3}{*}{\cite{WD1}} & {\textcolor{green}{\textbf{-}}  Studying weakly Douglas Finsler metric: $D_j{^i}_{kl}=T_{jkl}y^i$.}&\\
 && {} &  {\textcolor{green}{\textbf{-}} Randers metrics ($n\geq 3$): weakly Douglas  $\Leftrightarrow$ $D=0$.}&\\
 && {} &  {\textcolor{green}{\textbf{-}} Kropina metrics: weakly Douglas  $\Leftrightarrow$ $D=0$.}\\
 \cline{3-4}

 && \multirow{2}{*}{\cite{WD2}} & {\textcolor{green}{\textbf{-}} weakly Douglas warped product Finsler metric $\Leftrightarrow$ $D=0$.}&\\
 && {} &  {\textcolor{green}{\textbf{-}}weakly Douglas warped product Finsler metric:  $B=0$ $\Leftrightarrow$ $L=0$.  }\\
 \cline{3-4}

 && \multirow{2}{*}{\cite{WD3}} & {\textcolor{green}{\textbf{-}} under specific conditions, a Finsler warped product metrics: }&\\
 && {} &  {locally projectively flat $\Leftrightarrow$ it exhibits scalar flag Cur.  }\\
  \hline
\end{tabular}
\end{center}

\section{Conclusion and future directions}\label{conclusion}

The review paper is organized to provide a thorough overview of the historical development and significance of Douglas curvature in Finsler geometry. It is divided into several sections, each examining different aspects of Douglas curvature and its applications.
Each section distills the key findings and insights gained from the investigation of Douglas curvature in Finsler metrics and their generalizations, as presented in the tables. The accompanying tables present a concise summary of these findings, enabling researchers to quickly identify their areas of interest and navigate the research landscape related to Douglas curvature and its applications.

In the subsequent, we introduce several ideas for future research related to Douglas curvature in Finsler geometry. This section outlines promising directions for further exploration, emphasizing opportunities to deepen our understanding of this geometric concept and its potential applications.

- Proposing computational methods for the efficient calculation of Douglas curvature of different classes of Finsler metrics similar to the approach used for spherically symmetric in $\mathbb{R}^n$ and devising novel criteria and metrics for evaluating the quality and properties of Douglas curvature. We can draw inspiration from existing researches and adapt it to the specific context of Douglas curvature.

- Investigating the relationships between Douglas curvature and other curvature, that have not been extensively studied in the researches summarized in this paper.

- Exploring the behavior of Douglas curvature under various metric transformations, with a focus on Finsler metrics of Douglas type II, $GDW$-metrics, weakly Douglas metrics, and all other Finsler metrics related to Douglas metrics. Special attention should be given to transformations that have not been considered in the papers collected in Section \ref{Transformation} of this paper.

- Introducing new classes of Finsler metrics with interesting curvature properties through the study of Douglas curvature. The various generalizations of Douglas metrics proposed in different papers were collected in Section \ref{Generalization}. This provides a foundation for researchers to propose new classes of Finsler metrics that generalize Douglas metrics or GDW-metrics, which are projective invariant. By studying the properties and relationships of these new classes of metrics with the previously mentioned ones, we can further expand our understanding of Douglas curvature.

- Exploring the potential applications of Douglas metrics in physics. Notably, some papers have already been considered in this review as applications of Berwald spaces in Section \ref{Douglas}, such as \cite{FinsPhysics3} and \cite{GR1}. Given that Douglas spaces are a generalization of Berwald spaces, it is intriguing to consider whether these papers could also be relevant to Douglas spaces and their generalizations. By examining the connections between these papers and Douglas metrics or their generalizations, we may uncover new insights into the physical applications of Douglas spaces and their generalizations.

\end{document}